\documentclass[final,1p,times]{elsarticle}

\usepackage{amssymb}
 \usepackage{amsthm}
\usepackage{amscd}
\usepackage{amsmath}
\usepackage{amsfonts}
\usepackage{amssymb}
\usepackage{graphicx}
\newtheorem{theorem}{Theorem}
\newtheorem{proposition}[theorem]{Proposition}
\newtheorem{lemma}[theorem]{Lemma}
\newtheorem{corollary}[theorem]{Corollary}
\usepackage{mathrsfs}
\usepackage{titletoc}

\newcommand{\la}{\Delta}

\newcommand{\ra}{\rightarrow}
\newcommand{\p}{\partial}
\newcommand{\f}{\frac}

\newcommand{\be}{\begin{equation}}
\renewcommand{\ra}{\rightarrow}
\newcommand{\ee}{\end{equation}}
\newcommand{\bea}{\begin{eqnarray}}
\newcommand{\eea}{\end{eqnarray}}
\newcommand{\bna}{\begin{eqnarray*}}
\newcommand{\ena}{\end{eqnarray*}}

\renewcommand{\le}{\left}
\newcommand{\ri}{\right}

\journal{***}

\begin{document}

\begin{frontmatter}

\title{An improved Hardy-Trudinger-Moser inequality}

\author{Yunyan Yang}
 \ead{yunyanyang@ruc.edu.cn}
 \author{Xiaobao Zhu}
 \ead{zhuxiaobao@ruc.edu.cn}
\address{ Department of Mathematics,
Renmin University of China, Beijing 100872, P. R. China}

\begin{abstract}
Let $\mathbb{B}$ be the unit disc in $\mathbb{R}^2$, $\mathscr{H}$ be the completion of $C_0^\infty(\mathbb{B})$ under the norm
$$\|u\|_{\mathscr{H}}=\le(\int_\mathbb{B}|\nabla u|^2dx-\int_\mathbb{B}\frac{u^2}{(1-|x|^2)^2}dx\ri)^{1/2},\quad\forall u\in C_0^\infty(\mathbb{B}).$$
Denote $\lambda_1(\mathbb{B})=\inf_{u\in \mathscr{H},\,\|u\|_2=1}\|u\|_{\mathscr{H}}^2$,
where $\|\cdot\|_2$ stands for the $L^2(\mathbb{B})$-norm. Using blow-up analysis, we prove that for any $\alpha$,
$0\leq \alpha<\lambda_1(\mathbb{B})$,
$$\sup_{u\in\mathscr{H},\,\|u\|_{\mathscr{H}}^2-\alpha\|u\|_2^2\leq 1}\int_\mathbb{B} e^{4\pi u^2}dx<+\infty,$$
and that the above supremum can be attained by some function $u\in \mathscr{H}$ with
$\|u\|_{\mathscr{H}}^2-\alpha\|u\|_2^2= 1$.
This improves an earlier result of G. Wang and D. Ye \cite{Wang-Ye}.
  \end{abstract}

\begin{keyword}
Hardy-Trudinger-Moser inequality\sep
Trudinger-Moser inequality\sep blow-up analysis

\MSC[2010] 46E35
\end{keyword}

\end{frontmatter}

\titlecontents{section}[0mm]
                       {\vspace{.2\baselineskip}}
                       {\thecontentslabel~\hspace{.5em}}
                        {}
                        {\dotfill\contentspage[{\makebox[0pt][r]{\thecontentspage}}]}
\titlecontents{subsection}[3mm]
                       {\vspace{.2\baselineskip}}
                       {\thecontentslabel~\hspace{.5em}}
                        {}
                       {\dotfill\contentspage[{\makebox[0pt][r]{\thecontentspage}}]}

\setcounter{tocdepth}{2}


\section{Introduction}

Let $\mathbb{B}$ be the unit disc in $\mathbb{R}^2$ and $W_0^{1,2}(\mathbb{B})$ be the usual Sobolev
space. The Trudinger-Moser inequality \cite{24,19,17,22,Moser} says that for any $\beta\leq 4\pi$,
\be\label{Tr}\sup_{u\in W_0^{1,2}(\mathbb{B}),\,\|\nabla u\|_2\leq 1}\int_\mathbb{B} e^{\beta u^2}dx<\infty.\ee
Here and throughout this paper we denote the $L^p(\mathbb{B})$-norm by $\|\cdot\|_p$. This inequality is sharp in the sense that
for any $\beta>4\pi$, the integrals in (\ref{Tr}) are still finite but the supremum is infinity. It is a very powerful
tool in the problem of prescribed Gaussian curvature and partial differential equations.

Another important inequality in analysis is the Hardy inequality, namely
$$\label{Hardyineq}\int_{\mathbb{B}}|\nabla u|^2dx\geq\int_{\mathbb{B}}\f{u^2}{(1-|x|^2)^2}dx,\quad \forall u\in W_0^{1,2}(\mathbb{B}).$$
The Hardy inequality was improved in many ways. It was proved by H. Brezis and M. Marcus \cite{Brezis-Mar} that there exists some constant
$C$ such that
\be\label{Brz}\int_{\mathbb{B}}|\nabla u|^2dx-\int_{\mathbb{B}}\f{u^2}{(1-|x|^2)^2}dx\geq C\int_{\mathbb{B}} u^2dx,
\quad\forall u\in W_0^{1,2}(\mathbb{B}).\ee
Further improvements known as the Hardy-Sobolev inequalities were done by Maz'ya (\cite{Mazya}, Corollary 3, Section 2.1.6),
 Mancini-Sandeep \cite{Mancini-Sandeep}, Adimurthi-do \'O-Tintarev \cite{ADT}, and Mancini-Sandeep-Tintarev \cite{MST}.
In view of (\ref{Brz}),
$$\|u\|_{\mathscr{H}}=\le(
\int_{\mathbb{B}}|\nabla u|^2dx-\int_{\mathbb{B}}
\frac{u^2}{(1-|x|^2)^2}dx\ri)^{1/2}$$
defines a norm on $C_0^\infty(\mathbb{B})$. Let $\mathscr{H}$ be the completion of $C_0^\infty(\mathbb{B})$ under the norm
$\|\cdot\|_{\mathscr{H}}$. Clearly $\mathscr{H}$ is a Hilbert space.
By a result of Mancini-Sandeep (\cite{Mancini-Sandeep}, the inequality (1.2)), we can see that for any $p>1$,
there exists a constant $C_p>0$ such that
\be\label{2d-Hardy-Sobolev}
\|u\|_p\leq C_p\|u\|_{\mathscr{H}},\quad \forall u\in\mathscr{H}.
\ee
This is also obtained by Wang-Ye \cite{Wang-Ye}.
Thus we have
$$W_0^{1,2}(\mathbb{B})\subset\mathscr{H}\subset \cap_{p\geq 1}L^p(\mathbb{B}).$$
Obviously $\mathscr{H}\not\subset L^\infty(\mathbb{B})$. In view of (\ref{Tr}), one can expect
a Hardy-Trudinger-Moser inequality, namely
\be\label{htm}\sup_{u\in\mathscr{H},\,\|u\|_{\mathscr{H}}\leq 1}\int_{\mathbb{B}}e^{4\pi u^2}dx<+\infty.\ee
This was done by Wang-Ye by using blow-up analysis in \cite{Wang-Ye}, where the existence of extremal function
for (\ref{htm}) was also obtained. The inequality (\ref{htm}) was further extended by C. Tintarev \cite{Tint}
to a generalized Euclidean version by using Ground state transform, and by Mancini-Sandeep-Tintarev \cite{MST} to a hyperbolic
space version via a rearrangement argument. Compared with (\ref{htm}), another kind of singular Trudinger-Moser inequalities were obtained by Adimurthi-Sandeep
\cite{Adi-Sandeep}, Adimurthi-Yang \cite{Adi-Yang}, and de Souza-do \'O \cite{de-O}.

Motivated by the works of Adimurthi-Druet \cite{A-D}, Y. Yang \cite{Yang-JFA,Yang-Tran,Yang-MZ,Yang-JDE}
and C. Tintarev \cite{Tint},
we aim to rewrite (\ref{htm}) with $\|u\|_{\mathscr{H}}$ replaced by certain equivalent norm on $\mathscr{H}$.
To clarify this problem, we define
\be\label{eigenvalue}
\lambda_1(\mathbb{B})=\inf_{u\in\mathscr{H},\,u\not\equiv 0}\frac{\|u\|_\mathscr{H}^2}{\|u\|_2^2}.
\ee
By (\ref{2d-Hardy-Sobolev}) and a variational direct method, we have that $\lambda_1(\mathbb{B})$ can be attained by
some function $u\in\mathscr{H}$ with $\|u\|_{2}=1$. In particular, $\lambda_1(\mathbb{B})>0$.
In fact, $\lambda_1(\mathbb{B})$ is the first
eigenvalue of the Hardy-Laplace operator, namely
$$\mathscr{L}_H=-\Delta-\frac{1}{(1-|x|^2)^2}.$$
For any $\alpha$, $0\leq\alpha<\lambda_1(\mathbb{B})$ and any $u\in\mathscr{H}$, we denote
\be\label{norm-1}
\|u\|_{1,\alpha}=\le(\|u\|_{\mathscr{H}}^2-\alpha\|u\|_2^2\ri)^{1/2}.
\ee
Clearly $\|\cdot\|_{1,\alpha}$ is equivalent to $\|\cdot\|_{\mathscr{H}}$.
Our main result is the following:

\begin{theorem}\label{thm1}
Let $\mathbb{B}$ be the unit ball in $\mathbb{R}^2$, $\lambda_1(\mathbb{B})$ be defined as in (\ref{eigenvalue}).
Then for any $\beta\leq 4\pi$ and any $\alpha$, $0\leq\alpha<\lambda_1(\mathbb{B})$, the supremum
$$\sup_{u\in \mathscr{H},\,\|u\|_{1,\alpha}\leq 1}\int_\mathbb{B} e^{\beta u^2}dx$$
can be attained by some function $u_0\in \mathscr{H}$ with $\|u_0\|_{1,\alpha}=1$, where
$\|\cdot\|_{1,\alpha}$ is defined as in (\ref{norm-1}).
\end{theorem}

An interesting consequence of Theorem \ref{thm1} is the following weak form of the
Hardy-Trudinger-Moser inequality.

\begin{corollary}
Let $\lambda_1(\mathbb{B})$ be defined as in (\ref{eigenvalue}).
Then for any $\alpha$, $0\leq \alpha<\lambda_1(\mathbb{B})$,
there exists a constant $C>0$ depending only on $\alpha$ such that
$$\int_{\mathbb{B}}|\nabla u|^2dx-\int_\mathbb{B}\f{u^2}{(1-|x|^2)^2}dx
-\alpha\int_\mathbb{B}u^2dx-16\pi\log\int_\mathbb{B}e^udx\geq -C,\quad\forall u\in W_0^{1,2}(\mathbb{B}).$$
\end{corollary}

Following the lines of Y. Li \cite{Lijpde}, Adimurthi-Druet \cite{A-D}, Yang \cite{Yang-Tran},
and Wang-Ye \cite{Wang-Ye}, we prove Theorem \ref{thm1}
by using blow-up analysis.
We remark that Wang-Ye \cite{Wang-Ye} solved (\ref{htm}) by using
a result of Carleson-Chang \cite{CC} in addition to standard blow-up analysis. This method
was originally used by Li-Liu-Yang \cite{LLY}. In this paper, we
shall employ the capacity estimate introduced by Y. Li \cite{Lijpde}, instead of Carleson-Chang's result.
It would be interesting to extend our Theorem \ref{thm1} to the case involving $L^p(\mathbb{B})$-norm as in \cite{Lu-Yang}.

Earlier works in this direction were due to Carleson-Chang \cite{CC}, M. Struwe \cite{Struwe},
F. Flucher \cite{Flucher}, K. Lin \cite{Lin}, Ding-Jost-Li-Wang \cite{DJLW}, and Adimurthi-Struwe \cite{Adi-Stru}.
The remaining part of this paper is organized as follows. In Section 2, we give several preliminary lemmas; In Section 3,
we prove Theorem \ref{thm1}.

\section{Preliminary results}

In this section, we list several properties of the space $\mathscr{H}$.
Let
$$\mathscr{S}_0=\le\{u\in C_0^\infty(B): u(x)=u(r), u^\prime(r)\leq 0, \,{\rm where}\, r=|x|\ri\}$$
and $\mathscr{S}$ be the completion of $\mathscr{S}_0$ under the norm $\|\cdot\|_{\mathscr{H}}$. The following
embedding theorem was proved by Wang-Ye \cite{Wang-Ye}:

\begin{lemma}\label{L-2}
  $\mathscr{S}$ is embedded continuously in
$W_{\rm loc}^{1,2}(\mathbb{B})\cap C_{\rm loc}^{0,\f{1}{2}}(\mathbb{B}\setminus\{0\})$. Moreover,
$\mathscr{S}$ is embedded compactly in $L^p(\mathbb{B})$ for any $p\geq 1$.
\end{lemma}

The second important property of $\mathscr{H}$ is an embedding of Orlicz type, namely

\begin{lemma}\label{L-3}
 For any $p>1$ and any $u\in\mathscr{H}$, there holds
$$\int_{\mathbb{B}}e^{pu^2}dx<+\infty.$$
\end{lemma}

\noindent{\it Proof.} Fix $p>1$ and $u\in\mathscr{H}$. Since $C_0^\infty(\mathbb{B})$ is dense in
$\mathscr{H}$, we take $u_0\in C_0^\infty(\mathbb{B})$ such that
$\|u-u_0\|_{\mathscr{H}}^2<\pi/p$. Using an inequality $2ab\leq a^2+b^2$ twice, we have
\bna
\int_{\mathbb{B}}e^{pu^2}dx\leq \int_{\mathbb{B}}e^{4p(u-u_0)^2}dx+
\int_{\mathbb{B}}e^{4pu_0^2}dx.
\ena
By (\cite{Wang-Ye}, Theorem 1), we have
$$\int_{\mathbb{B}}e^{4p(u-u_0)^2}dx<+\infty.$$
Since $u_0$ is uniformly bounded in $\mathbb{B}$, we have
$\int_{\mathbb{B}}e^{4pu_0^2}dx<+\infty$.
This gives the desired estimate. $\hfill\Box$\\

Finally we state an obvious but very important property of $\mathscr{H}$.

\begin{lemma}\label{L-4}
 Suppose $u\in W^{1,2}_{\rm loc}(\mathbb{B})$, $v\in\mathscr{H}$ and
$v$ is radially symmetric.
If there exists some $r$, $0<r<1$, such that $u=v$ on $\mathbb{B}\setminus\mathbb{B}_r$, then
$u\in\mathscr{H}$.
\end{lemma}

\noindent{\it Proof.} It follows from Lemma \ref{L-2} that $v\in W_{\rm loc}^{1,2}(\mathbb{B})$. Clearly we have
$u-v\in W_0^{1,2}(\mathbb{B}_r)\subset W_0^{1,2}(\mathbb{B})$. This leads to $u=(u-v)+v\in \mathscr{H}$.
$\hfill\Box$

\section{Proof of Theorem \ref{thm1}}

In this section, we prove Theorem 1 by using a blow-up scheme similar to that of Wang-Ye \cite{Wang-Ye}, and
thereby follow closely Y. Li \cite{Lijpde}, Adimurthi-Druet \cite{A-D}, and Yang \cite{Yang-Tran}.
The proof will be divided into several subsections.

\subsection{The subcritical case\\}

In this subsection, we prove that the subcritical Trudinger-Moser functional $J(u)=\int_{\mathbb{B}} e^{\gamma u^2}dx$
has a maximizer in the function space $\{u\in \mathscr{H}: \|u\|_{1,\alpha}\leq 1\}$  for
any $\gamma<4\pi$ and any $\alpha$, $0\leq\alpha<\lambda_1(\mathbb{B})$.

\begin{proposition}\label{P-1}
 Let $0\leq\alpha<\lambda_1(\mathbb{B})$.
For any $\epsilon>0$, there exists  $u_\epsilon\in\mathscr{S}\cap C^\infty(\mathbb{B})
\cap
C^0(\overline{\mathbb{B}})$ such that
$\|u_\epsilon\|_{1,\alpha}=1$ and
\be\label{inequa}\int_{\mathbb{B}}e^{(4\pi-\epsilon)u_\epsilon^2}dx=\sup_{u\in\mathscr{H},\,\|u\|_{1,\alpha}\leq 1}
\int_{\mathbb{B}}e^{(4\pi-\epsilon)u^2}dx.\ee
Moreover, $u_\epsilon$ satisfies the following Euler-Lagrange equation
\be\label{Euler-Lagrange}
\le\{\begin{array}{lll}
-\Delta u_\epsilon-\frac{u_\epsilon}{(1-|x|^2)^2}-\alpha u_\epsilon=\f{1}{\lambda_\epsilon}u_\epsilon e^{(4\pi-\epsilon)u_\epsilon^2}
\,\,{\rm in}\,\, \mathbb{B},\\[1.2ex]
\lambda_\epsilon=\int_{\mathbb{B}}u_\epsilon^2e^{(4\pi-\epsilon)u_\epsilon^2}dx.
\end{array}
\ri.
\ee
Furthermore, we have
\be\label{lambda-epsilon}\liminf_{\epsilon\ra 0}\lambda_\epsilon>0.\ee
\end{proposition}
\noindent{\it Proof}.
 We first claim that for any $\gamma<4\pi$, there holds
\be\label{H-inequa}
\sup_{u\in\mathscr{H},\,\|u\|_{\mathscr{H}}\leq 1}\int_{\mathbb{B}}e^{\gamma u^2}dx
<+\infty.
\ee
To see this, we use an argument of radially decreasing rearrangement with respect to the standard hyperbolic metric $dv=\f{1}{(1-|x|^2)^2}dx$.
For any $u\in C_0^\infty(\mathbb{B})$, we let $u^\ast$ be the radially decreasing rearrangement of $|u|$ with respect to the standard hyperbolic metric.
It follows
from \cite{Baernstein} that
\bna
\int_{\mathbb{B}}|\nabla u^\ast|^2dx\leq\int_{\mathbb{B}}|\nabla u|^2dx,\quad\int_{\mathbb{B}}\frac{{u^\ast}^2}{(1-|x|^2)^2}dx=\int_{\mathbb{B}}\frac{{u}^2}{(1-|x|^2)^2}dx.
\ena
Clearly we have
$\|u^\ast\|_{\mathscr{H}}\leq \|u\|_{\mathscr{H}}$. This leads to
$$\sup_{u\in C_0^\infty(\mathbb{B}),\,\|u\|_{\mathscr{H}}\leq 1}\int_{\mathbb{B}}e^{\gamma u^2}dx
\leq
\sup_{u\in\mathscr{S},\,\|u\|_{\mathscr{H}}\leq 1}\int_{\mathbb{B}}e^{\gamma u^2}dx.$$
In view of Lemma \ref{L-2}, for any $u\in\mathscr{H}$ with $\|u\|_{\mathscr{H}}\leq 1$, there exists a sequence of functions
$u_j\in C_0^\infty(\mathbb{B})$ such that $\|u_j\|_{\mathscr{H}}\leq 1$, $u_j\ra u$ in $\mathscr{H}$ and $u_j\ra u$ a. e.
in $\mathbb{B}$. Thus
$$\int_{\mathbb{B}}e^{\gamma u^2}dx\leq\limsup_{j\ra\infty}\int_{\mathbb{B}}e^{\gamma u_j^2}dx.$$
Hence
$$\sup_{u\in\mathscr{H},\,\|u\|_{\mathscr{H}}\leq 1}\int_{\mathbb{B}}e^{\gamma u^2}dx=
\sup_{u\in\mathscr{S},\,\|u\|_{\mathscr{H}}\leq 1}\int_{\mathbb{B}}e^{\gamma u^2}dx,$$
which together with (\cite{Wang-Ye}, Theorem 3) implies (\ref{H-inequa}).

Now let $0<\epsilon<4\pi$ be fixed. Note that
\bna
\int_{\mathbb{B}}{u^\ast}^2dx&=&\int_{\mathbb{B}}{u^\ast}^2(1-|x|^2)^2dv\\
&\geq&\int_{\mathbb{B}}(u^2(1-|x|^2)^2)^\ast dv\\
&=&\int_{\mathbb{B}}u^2(1-|x|^2)^2dv\\
&=&\int_{\mathbb{B}}{u}^2dx,
\ena
and that
\bna
\int_{\mathbb{B}}e^{(4\pi-\epsilon){u^\ast}^2}dx&=&\int_{\mathbb{B}}e^{(4\pi-\epsilon){u^\ast}^2}(1-|x|^2)^2dv\\
&\geq&\int_{\mathbb{B}}e^{(4\pi-\epsilon){u}^2}(1-|x|^2)^2dv\\
&=&\int_{\mathbb{B}}e^{(4\pi-\epsilon){u}^2}dx,
\ena
where  we have used the Hardy-Littlewood inequality (see \cite{Brock}) and the fact that
the radially decreasing rearrangement of $(1-|x|^2)^2$ with respect to the standard hyperbolic metric is itself. Therefore
$$\sup_{u\in C_0^\infty(\mathbb{B}),\,\|u\|_{1,\alpha}\leq 1}
\int_{\mathbb{B}}e^{(4\pi-\epsilon)u^2}dx\leq\sup_{u\in\mathscr{S},\,\|u\|_{1,\alpha}\leq 1}
\int_{\mathbb{B}}e^{(4\pi-\epsilon)u^2}dx.$$
Since $C_0^\infty(\mathbb{B})$ is dense in $\mathscr{H}$, we obtain
\be\label{subcritical}\sup_{u\in \mathscr{H},\,\|u\|_{1,\alpha}\leq 1}
\int_{\mathbb{B}}e^{(4\pi-\epsilon)u^2}dx=\sup_{u\in\mathscr{S},\,\|u\|_{1,\alpha}\leq 1}
\int_{\mathbb{B}}e^{(4\pi-\epsilon)u^2}dx.\ee

To prove (\ref{inequa}), we use a method of variation. Observing (\ref{H-inequa}) and (\ref{subcritical}), we can take a sequence of
functions $u_j\in\mathscr{S}$ with $\|u_j\|_{1,\alpha}\leq 1$ such that
\be\label{2alp}
\int_{\mathbb{B}}e^{(4\pi-\epsilon)u_j^2}dx\ra\sup_{u\in \mathscr{H},\,\|u\|_{1,\alpha}\leq 1}
\int_{\mathbb{B}}e^{(4\pi-\epsilon)u^2}dx\quad{\rm as}\quad j\ra \infty.\ee
Since $0\leq\alpha<\lambda_1(\mathbb{B})$, we have
$$\|u_j\|_{\mathscr{H}}^2\leq \frac{\lambda_1(\mathbb{B})}{\lambda_1(\mathbb{B})-\alpha}.$$
Note that $\mathscr{H}$ is a Hilbert space. Up to a subsequence there exists some $u_\epsilon\in\mathscr{S}$ such that
\bna
&u_j\rightharpoonup u_\epsilon&{\rm weakly\,\,in}\quad \mathscr{H},\\
&u_j\ra u_\epsilon &{\rm strongly\,\,in}\quad L^p(\mathbb{B}),\,\,\forall p\geq 1,\\
&u_j\ra u_\epsilon &{\rm a.\, e.\,\,in} \quad\mathbb{B}.
\ena
It follows that
\be\label{leq-1}\|u_\epsilon\|_{1,\alpha}^2\leq\liminf_{j\ra\infty}\|u_j\|_{1,\alpha}^2\leq 1.\ee
A straightforward calculation shows
\bea
\|u_j-u_\epsilon\|_{\mathscr{H}}^2&=&\langle u_j-u_\epsilon,u_j-u_\epsilon \rangle_{\mathscr{H}}{\nonumber}\\
&=&\|u_j\|_{\mathscr{H}}^2+\|u_\epsilon\|_{\mathscr{H}}^2-2\langle u_j,u_\epsilon \rangle_{\mathscr{H}}{\nonumber}\\
&=&\|u_j\|_{\mathscr{H}}^2-\|u_\epsilon\|_{\mathscr{H}}^2+o_j(1){\nonumber}\\
&=&\|u_j\|_{1,\alpha}^2-\|u_\epsilon\|_{1,\alpha}^2+o_j(1),\label{lt-1}
\eea
since $\langle u_j,u_\epsilon \rangle_{\mathscr{H}}\ra \|u_\epsilon\|_{\mathscr{H}}^2$ and $\|u_j\|_2\ra
\|u_\epsilon\|_2$ as $j\ra \infty$.

For any $\nu>0$, using an elementary inequality $2ab\leq \nu a^2+b^2/\nu$, we have
$$\label{shou}u_j^2\leq (1+\nu)(u_j-u_\epsilon)^2+(1+{1}/{\nu})u_\epsilon^2.$$
Choosing $\nu=\epsilon/(8\pi-2\epsilon)$ in the above equation, we have
\be\label{d-1}(4\pi-\epsilon)u_j^2\leq (4\pi-\epsilon/2)(u_j-u_\epsilon)^2+\frac{32\pi^2}{\epsilon}u_\epsilon^2.\ee
From (\ref{lt-1}) we can find some positive integer $j_0$ such that
$$\|u_j-u_\epsilon\|_{\mathscr{H}}^2\leq \f{4\pi-\epsilon/3}{4\pi-\epsilon/2},\quad\forall j\geq j_0.$$
This together with (\ref{d-1}) gives
\be\label{inequa-0}(4\pi-\epsilon)u_j^2\leq (4\pi-\epsilon/3)\frac{(u_j-u_\epsilon)^2}{\|u_j-u_\epsilon\|_{\mathscr{H}}^2}
+\frac{32\pi^2}{\epsilon}u_\epsilon^2,\quad\forall j\geq j_0.\ee
By Lemma \ref{L-3},
\be\label{in-1}\int_{\mathbb{B}}e^{q u_\epsilon^2}dx<\infty,\quad \forall q>1.\ee
Take
$$p=\f{4\pi-\epsilon/4}{4\pi-\epsilon/3}>1.$$
Combining (\ref{H-inequa}), (\ref{inequa-0}) and (\ref{in-1}), we conclude that $e^{(4\pi-\epsilon)u_j^2}$ is bounded in
$L^p(\mathbb{B})$, which together with $u_j\ra u_0$ a. e. as $j\ra\infty$ implies that
$e^{(4\pi-\epsilon)u_j^2}$ converges to $e^{(4\pi-\epsilon)u_\epsilon^2}$  in $L^1(\mathbb{B})$.
This together with (\ref{2alp}) leads to (\ref{inequa}).
Recall (\ref{leq-1}) we have $\|u_\epsilon\|_{1,\alpha}\leq 1$. Now we claim that $\|u_\epsilon\|_{1,\alpha}=1$. For otherwise, we have
$\|u_\epsilon\|_{1,\alpha}<1$, and thus
\bna
\sup_{u\in\mathscr{H},\,\|u\|_{1,\alpha}\leq 1}
\int_{\mathbb{B}}e^{(4\pi-\epsilon)u^2}dx&=&\int_{\mathbb{B}}e^{(4\pi-\epsilon)u_\epsilon^2}dx\\
&<&\int_{\mathbb{B}}e^{(4\pi-\epsilon){u_\epsilon^2}/
{\|u_\epsilon\|_{1,\alpha}^2}}dx\\&\leq&\sup_{u\in\mathscr{H},\,\|u\|_{1,\alpha}\leq 1}
\int_{\mathbb{B}}e^{(4\pi-\epsilon)u^2}dx,\ena
which is a contradiction.

It is not difficult to see that $u_\epsilon$ satisfies the Euler-Lagrange equation (\ref{Euler-Lagrange}).
Finally $u_\epsilon\in C^\infty(\mathbb{B})$ follows from standard elliptic estimates,
and the fact that $u_\epsilon\in C^0(\overline{\mathbb{B}})$ follows from $u_\epsilon\in\mathscr{S}$.

Finally we prove (\ref{lambda-epsilon}). Using an elementary inequality $e^t\leq 1+te^t$ for $t\geq 0$, we have
$$\int_{\mathbb{B}}e^{(4\pi-\epsilon)u_\epsilon^2}dx\leq \pi+(4\pi-\epsilon)\lambda_\epsilon.
$$
Note that $\int_{\mathbb{B}}e^{(4\pi-\epsilon)u^2}dx$ is monotone with respect to $\epsilon>0$.
 For any fixed $u\in\mathscr{H}$ with $\|u\|_{1,\alpha}=1$, in view of Lemma \ref{L-3}, there holds
$$\pi<\int_{\mathbb{B}}e^{4\pi u^2}dx=\lim_{\epsilon\ra 0}\int_{\mathbb{B}}e^{(4\pi-\epsilon)u^2}dx
\leq\liminf_{\epsilon\ra 0}\int_{\mathbb{B}}e^{(4\pi-\epsilon)u_\epsilon^2}dx\leq \pi+4\pi\liminf_{\epsilon\ra 0}
\lambda_\epsilon.$$
This leads to (\ref{lambda-epsilon}) immediately.
$\hfill\Box$

\subsection{Blow-up analysis\\}

In this subsection, we perform the blow-up procedure.
Let $u_\epsilon$ be as in Proposition \ref{P-1}.
Since $\|u_\epsilon\|_{1,\alpha}=1$ and $\alpha<\lambda_1(\mathbb{B})$,
$u_\epsilon$ is bounded in $\mathscr{H}$. By Lemma \ref{L-2},  there exists $u_0\in L^2(\mathbb{B})$ such that
 up to a
 subsequence,
$u_\epsilon\ra u_0$  in $L^2(\mathbb{B})$
and $u_\epsilon\ra u_0$ a. e. in $\mathbb{B}$  as
$\epsilon\ra 0$.
On the other hand, there exists some $v_0\in\mathscr{S}$ such that $u_\epsilon\rightharpoonup v_0$
weakly in $\mathscr{H}$. In particular
 $$\lim_{\epsilon\ra 0}\int_{\mathbb{B}}u_\epsilon\varphi dx=\int_{\mathbb{B}}v_0\varphi dx,\quad\forall
 \varphi\in L^2(\mathbb{B}).$$
 Since the weak limit is unique, we have $v_0=u_0$. In conclusion, there exists $u_0\in\mathscr{S}$
 such that up to a subsequence,
 $$u_\epsilon\rightharpoonup u_0\,\,\,{\rm in}\,\,\, \mathscr{H},\,\,\,
 u_\epsilon\ra u_0\,\,\,{\rm a.\,e.\,\,in}\,\,\,\mathbb{B}$$
 as $\epsilon\ra 0$. Noting that $\langle u_\epsilon,u_0\rangle_{\mathscr{H}}\ra \langle u_0,u_0\rangle_{\mathscr{H}}$,
 we have
 $$\|u_0\|_{1,\alpha}^2=\|u_0\|_{\mathscr{H}}^2-\alpha\|u_0\|_2^2\leq \liminf_{\epsilon\ra 0}
 \|u_\epsilon\|_{1,\alpha}=1.$$

Let $c_\epsilon=u_\epsilon(0)=\max_{\mathbb{B}}u_\epsilon$.
If $c_\epsilon$ is bounded, we have by using the Lebesgue dominated convergence theorem,
$$\int_\mathbb{B}e^{4\pi u_0^2}dx=
\lim_{\epsilon\ra 0}\int_\mathbb{B}e^{(4\pi-\epsilon) u_\epsilon^2}dx=
\sup_{u\in\mathscr{H},\, \|u\|_{1,\alpha}\leq 1}\int_\mathbb{B}e^{4\pi u^2}dx.$$
Hence $u_0$ is a desired extremal function and Theorem 1 holds. In the following, we assume
\be\label{cepsil}c_\epsilon\ra +\infty\,\,{\rm as}\,\,\epsilon\ra 0.\ee
Now we claim that $u_0\equiv0$. To see this, suppose $u_0\not\equiv 0$, then $\|u_0\|_{1,\alpha}>0$. On one hand, by the H\"older
inequality, $\forall \nu>0$, there holds
\be\label{H-d}
 \|e^{(4\pi-\epsilon)u_\epsilon^2}\|_{1+\nu}\leq
 \|e^{(4\pi-\epsilon)(u_\epsilon-u_0)^2}\|_{(1+\nu)(1+2\nu)}^{1+\nu}
 \|e^{(4\pi-\epsilon)u_0^2}\|_{(1+\nu)^2(1+2\nu)/\nu^2}^{1+1/\nu}.
\ee
On the other hand, we calculate
\bea
\|u_\epsilon-u_0\|_{\mathscr{H}}^2&=&\langle u_\epsilon-u_0,u_\epsilon-u_0\rangle_{\mathscr{H}}\nonumber\\
&=&\|u_\epsilon\|_{\mathscr{H}}^2-
\|u_0\|_{\mathscr{H}}^2+o_\epsilon(1)\nonumber\\
&=&\|u_\epsilon\|_{1,\alpha}^2-\|u_0\|_{1,\alpha}^2+o_\epsilon(1)\nonumber\\
&=& 1-\|u_0\|_{1,\alpha}^2+o_\epsilon(1)\nonumber\\
&<&1-\|u_0\|_{1,\alpha}^2/2,\label{Lio}
\eea
provided that $\epsilon$ is sufficiently small.
Choosing  $\nu=\|u_0\|_{1,\alpha}^2/16$ in (\ref{H-d}), we have
by (\ref{Lio}) and (\ref{H-inequa}) that $e^{(4\pi-\epsilon)u_\epsilon^2}$ is bounded in $L^{1+\nu}(\mathbb{B})$.
Then applying standard elliptic estimates to (\ref{Euler-Lagrange}), we get that $u_\epsilon$ is
 bounded in $C^0_{\rm loc}(\mathbb{B})$, which contradicts (\ref{cepsil}). Therefore $u_0\equiv 0$.

We set
$$\label{r-eps} r_\epsilon=\sqrt{\lambda_\epsilon}c_\epsilon^{-1}e^{-(2\pi-\epsilon/2)c_\epsilon^2}.$$
 For any $0<\delta<4\pi$, we have by using the H\"older inequality and (\ref{H-inequa}),
 $$\label{l-eps}\lambda_\epsilon=\int_{\mathbb{B}} u_\epsilon^2 e^{(4\pi-\epsilon)u_\epsilon^2}dx\leq e^{\delta c_\epsilon^2}
 \int_{\mathbb{B}} u_\epsilon^2 e^{(4\pi-\epsilon-\delta)u_\epsilon^2}dx\leq Ce^{\delta c_\epsilon^2}$$
 for some constant $C$ depending only on $\delta$.
 This leads to
 \be\label{r-0} r_\epsilon^2\leq C c_\epsilon^{-2}e^{-(4\pi-\epsilon-\delta)c_\epsilon^2}\ra 0\quad{\rm as}\quad \epsilon\ra 0.\ee
 Define two blow-up sequences of functions on $\mathbb{B}_{r_\epsilon^{-1}}=\{x\in\mathbb{R}^2: |x|< r_\epsilon^{-1}\}$ as
 $$
   \psi_\epsilon(x)=c_\epsilon^{-1}u_\epsilon(r_\epsilon
   x),\quad
   \varphi_\epsilon(x)=c_\epsilon(u_\epsilon(r_\epsilon
   x)-c_\epsilon).
   $$
   This kind of blow-up functions are suitable for such a problem was first discovered by Adimurthi-Struwe \cite{Adi-Stru}.
   A direct computation shows
   \be\label{p-s-eq}
   -\la \psi_\epsilon=\frac{r_\epsilon^2}{(1-r_\epsilon^2|x|^2)^2}\psi_\epsilon+\alpha r_\epsilon^2\psi_\epsilon+c_\epsilon^{-2}\psi_\epsilon e^{(4\pi-\epsilon)(1+\psi_\epsilon)\varphi_\epsilon}\quad{\rm
   in}\quad \mathbb{B}_{r_\epsilon^{-1}},
   \ee
   \be\label{phi-eq}
   -\la\varphi_\epsilon=\frac{r_\epsilon^2c_\epsilon^2}{(1-r_\epsilon^2|x|^2)^2} \psi_\epsilon+\alpha r_\epsilon^2c_\epsilon^2\psi_\epsilon
   +\psi_\epsilon e^{(4\pi-\epsilon)(1+\psi_\epsilon)\varphi_\epsilon}\quad{\rm
   in}\quad \mathbb{B}_{r_\epsilon^{-1}}.
   \ee
   We now consider the asymptotic behavior of $\psi_\epsilon$ and $\varphi_\epsilon$.
   By (\ref{r-0}), we have $r_\epsilon^2c_\epsilon^q\ra 0$ as $\epsilon\ra 0$ for any $q\geq 1$.
   Since $\mathbb{B}_{r_\epsilon^{-1}}\ra \mathbb{R}^2$ as $\epsilon\ra 0$, we have that $|\psi_\epsilon|\leq 1$ and
   $\Delta\psi_\epsilon(x)\ra 0$ uniformly in $\mathbb{B}_{R}$ for any fixed $R>0$ as $\epsilon\ra 0$. Applying elliptic estimates
   to (\ref{p-s-eq}), we have
   $\psi_\epsilon\ra \psi$ in $C^1_{\rm loc}(\mathbb{R}^2)$, where $\psi$ is a bounded harmonic function in $\mathbb{R}^2$.
   Note that $\psi(0)=\lim_{\epsilon\ra 0}\psi_\epsilon(0)=1$. The Liouville theorem implies that $\psi\equiv 1$ on
   $\mathbb{R}^2$. Hence
   $$\psi_\epsilon\ra 1\quad{\rm in}\quad C^1_{\rm loc}(\mathbb{R}^2).$$
    Since  $\varphi_\epsilon(x)\leq \varphi_\epsilon(0)=0$ for all $x\in \mathbb{B}_{r_\epsilon^{-1}}$,
    it is not difficult to see that $\Delta \varphi_\epsilon$ is uniformly bounded in $\mathbb{B}_{R}$
    for any fixed $R>0$.
    We then conclude by applying elliptic estimates to
   the equation (\ref{phi-eq}) that
   \be\label{phibubble}
   \varphi_\epsilon\ra \varphi\quad{\rm in}\quad
   C^1_{\rm{loc}}(\mathbb{R}^2),
   \ee
   where $\varphi$ satisfies
   $$\label{bubbel}
   \le\{\begin{array}{lll}&\la
     \varphi=-e^{8\pi\varphi}\quad{\rm in}\quad\mathbb{R}^2\\[1.2 ex]
     &\varphi(0)=0=\sup_{\mathbb{R}^2}\varphi\\[1.2 ex] &\int_{\mathbb{R}^2}e^{8\pi\varphi}dx\leq 1.
     \end{array}
    \ri.$$
    By a result of Chen-Li \cite{CL}, we have
    \be\label{v-r}\varphi(x)=-\f{1}{4\pi}\log(1+\pi|x|^2),\,\,\,\int_{\mathbb{R}^2}e^{8\pi\varphi}dx=1.\ee

   Now we consider the convergence behavior of $u_\epsilon$ away from zero.
   Set $u_{\epsilon,\beta}=\min\{u_\epsilon,\beta c_\epsilon\}$ for any $\beta$, $0<\beta<1$. Then we have

   \begin{lemma}\label{L-5}
    For any $\beta$, $0<\beta<1$, there holds $\lim_{\epsilon\ra 0}
   \|u_{\epsilon,\beta}\|_{1,\alpha}^2=\beta$.
   \end{lemma}

  \noindent{\it Proof}. Note that $(u_\epsilon-\beta c_\epsilon)^+\in W_0^{1,2}(\mathbb{B})$ and thus
  $u_{\epsilon,\beta}=u_\epsilon-(u_\epsilon-\beta c_\epsilon)^+\in \mathscr{H}$. Testing the equation (\ref{Euler-Lagrange})
  by $u_{\epsilon,\beta}$, we have
  $$\int_{\mathbb{B}}\le(\nabla u_{\epsilon,\beta}\nabla u_\epsilon-\frac{u_{\epsilon,\beta}u_\epsilon}{(1-|x|^2)^2}-\alpha
  u_{\epsilon,\beta}u_\epsilon\ri)dx=\int_{\mathbb{B}}\f{1}{\lambda_\epsilon}u_{\epsilon,\beta}u_\epsilon e^{(4\pi-\epsilon)u_\epsilon^2}dx.$$
  It follows that
  \bna
   \|u_{\epsilon,\beta}\|_{1,\alpha}^2&=&\int_{\mathbb{B}}\le(\nabla u_{\epsilon,\beta}\nabla u_\epsilon-\frac{u_{\epsilon,\beta}u_\epsilon}{(1-|x|^2)^2}-\alpha
  u_{\epsilon,\beta}u_\epsilon\ri)dx\\&&+\int_{\mathbb{B}}\f{u_{\epsilon,\beta}(u_\epsilon-u_{\epsilon,\beta})}{(1-|x|^2)^2}dx+
  \alpha\int_{\mathbb{B}}u_{\epsilon,\beta}(u_\epsilon-u_{\epsilon,\beta})dx\\
  &\geq&\int_{\mathbb{B}}\f{1}{\lambda_\epsilon}u_{\epsilon,\beta}u_\epsilon e^{(4\pi-\epsilon)u_\epsilon^2}dx\\
  &\geq&\beta\int_{\mathbb{B}_{R}(0)}(1+o_\epsilon(1))e^{8\pi\varphi}dy\ena
  for any $R>0$. Letting $\epsilon\ra 0$ first, and then $R\ra\infty$, we obtain
  $$\liminf_{\epsilon\ra 0}\|u_{\epsilon,\beta}\|_{1,\alpha}^2\geq \beta.$$

  Similarly, testing the equation (\ref{Euler-Lagrange}) by $(u_\epsilon-\beta c_\epsilon)^+$, we have
  \bna
  \|(u_\epsilon-\beta c_\epsilon)^+\|_{1,\alpha}^2&=&\int_{\mathbb{B}}\le(\nabla(u_\epsilon-\beta c_\epsilon)^+\nabla u_\epsilon-
  \frac{(u_\epsilon-\beta c_\epsilon)^+u_\epsilon}{(1-|x|^2)^2}-\alpha
  (u_\epsilon-\beta c_\epsilon)^+u_\epsilon\ri)dx\\&&+\int_{\mathbb{B}}\f{(u_\epsilon-\beta c_\epsilon)^+u_{\epsilon,\beta}}{(1-|x|^2)^2}dx+
  \alpha\int_{\mathbb{B}}(u_\epsilon-\beta c_\epsilon)^+u_{\epsilon,\beta}dx\\
  &\geq&\int_{\mathbb{B}}\f{1}{\lambda_\epsilon}(u_\epsilon-\beta c_\epsilon)^+u_\epsilon e^{(4\pi-\epsilon)u_\epsilon^2}dx\\
  &\geq&(1-\beta)\int_{\mathbb{B}_{R}(0)}(1+o_\epsilon(1))e^{8\pi\varphi}dy.
  \ena
  This implies that
  $$\liminf_{\epsilon\ra 0}\|(u_\epsilon-\beta c_\epsilon)^+\|_{1,\alpha}^2\geq (1-\beta).$$
  Since $u_\epsilon\ra 0$  in $L^p(\mathbb{B})$ as $\epsilon\ra 0$ for any fixed $p>1$, one can see that
  $$\lim_{\epsilon\ra 0}\le(\|u_{\epsilon,\beta}\|_{1,\alpha}^2+\|(u_\epsilon-\beta c_\epsilon)^+\|_{1,\alpha}^2-
  \|u_\epsilon\|_{1,\alpha}^2\ri)=0.$$
  Therefore
  $$
  \lim_{\epsilon\ra 0}
   \|u_{\epsilon,\beta}\|_{1,\alpha}^2=\beta,\,\,\,\lim_{\epsilon\ra 0}
   \|(u_\epsilon-\beta c_\epsilon)^+\|_{1,\alpha}^2=1-\beta.$$
   This completes the proof of the lemma. $\hfill\Box$\\

  \begin{lemma}\label{L-6}
    There holds
   $$\lim_{\epsilon\ra 0} \int_\mathbb{B}e^{(4\pi-\epsilon)u_\epsilon^2}dx=\pi+\limsup_{\epsilon\ra 0}
   \f{\lambda_\epsilon}{c_\epsilon^2}.$$
   \end{lemma}

  \noindent{\it Proof.} On one hand we have for any $\beta$, $0<\beta<1$,
  \bna
  \int_\mathbb{B}e^{(4\pi-\epsilon)u_\epsilon^2}dx&=&\int_{u_\epsilon<\beta
  c_\epsilon}e^{(4\pi-\epsilon)u_\epsilon^2}dx+
  \int_{u_\epsilon\geq\beta
  c_\epsilon}e^{(4\pi-\epsilon)u_\epsilon^2}dx\\
  &\leq&
  \int_\mathbb{B} e^{(4\pi-\epsilon)u_{\epsilon,\beta}^2}dx
  +\f{\lambda_\epsilon}{\beta^2 c_\epsilon^{2}}.
  \ena
  It follows from Lemma \ref{L-5} that
  $\int_\mathbb{B} e^{(4\pi-\epsilon)u_{\epsilon,\beta}^2}dx
  \ra |\mathbb{B}|=\pi$
  as $\epsilon\ra 0$. Hence
  $$\int_\mathbb{B}e^{(4\pi-\epsilon)u_\epsilon^2}dx\leq \pi+
  \f{\lambda_\epsilon}{\beta^2 c_\epsilon^{2}}+o_\epsilon(1).$$
  Letting $\epsilon\ra 0$ first, then $\beta\ra 1$ in the above inequality,
   we get
   \be\label{left}\lim_{\epsilon\ra 0}\int_\mathbb{B}e^{(4\pi-\epsilon)u_\epsilon^2}dx\leq \pi+
   \limsup_{\epsilon\ra 0}\f{\lambda_\epsilon}{c_\epsilon^2}.\ee

   On the other hand we have by (\ref{phibubble})
   $$\int_{\mathbb{B}_{Rr_\epsilon}}e^{(4\pi-\epsilon)u_\epsilon^2}dx=\f{\lambda_\epsilon}{c_\epsilon^2}
   \le(\int_{\mathbb{B}_R}e^{8\pi\varphi}dx+o_\epsilon(1)\ri).$$
   It is easy to see that
   $$\int_{\mathbb{B}_{Rr_\epsilon}}e^{(4\pi-\epsilon)u_\epsilon^2}dx
   \leq \int_{\mathbb{B}}e^{(4\pi-\epsilon)u_\epsilon^2}dx-\pi(1-R^2r_\epsilon^2).$$
   Combining the above two estimates and letting $\epsilon\ra 0$ first, then $R\ra +\infty$, we have
   \be\label{right}
   \limsup_{\epsilon\ra 0}\f{\lambda_\epsilon}{c_\epsilon^2}
   \leq\lim_{\epsilon\ra 0}\int_{\mathbb{B}}e^{(4\pi-\epsilon)u_\epsilon^2}dx-\pi.
   \ee
   Combining (\ref{left}) and (\ref{right}), we get the desired result. $\hfill\Box$\\

   Obviously Lemma \ref{L-6} implies that
   \be\label{t-z}\lim_{\epsilon\ra
  0}c_\epsilon/\lambda_\epsilon=0.\ee
  (Here and in the sequel we do {\it not}
distinguish sequence and subsequence.)
  This will be used to prove the following:

  \begin{lemma}\label{L-7} $\forall \phi\in
  C^\infty(\overline{\mathbb{B}})$,  we have
  $$\lim_{\epsilon\ra 0}\int_\mathbb{B} \phi\f{1}{\lambda_\epsilon}
  c_\epsilon u_\epsilon
  e^{(4\pi-\epsilon)u_\epsilon^2}dx=\phi(0).$$
  \end{lemma}

  \noindent{\it Proof.} For any fixed $\beta$, $0<\beta<1$, we divide $\mathbb{B}$ into three parts
  $$\mathbb{B}=(\{u_\epsilon>\beta c_\epsilon\}\setminus \mathbb{B}_{Rr_\epsilon})\cup
  (\{u_\epsilon\leq{\beta c_\epsilon}\}\setminus\mathbb{B}_{Rr_\epsilon})
  \cup \mathbb{B}_{Rr_\epsilon}.$$ Denote the integrals on the above three
  domains by $I_1$, $I_2$ and $I_3$ respectively. Firstly we have
  \bna
  |I_1|&\leq&\sup_\mathbb{B}|\phi|\int_{\{u_\epsilon>{\beta c_\epsilon}\}
  \setminus \mathbb{B}_{Rr_\epsilon}}\f{1}{\lambda_\epsilon}c_\epsilon
  u_\epsilon e^{(4\pi-\epsilon)u_\epsilon^2}dx\\
  &\leq&\f{1}{\beta}\sup_\mathbb{B}|\phi|\le(1-\int_{\mathbb{B}_{Rr_\epsilon}}
  \f{1}{\lambda_\epsilon}
  u_\epsilon^2e^{(4\pi-\epsilon)u_\epsilon^2}dx\ri)\\
  &=&\f{1}{\beta}\sup_\mathbb{B}|\phi|
  \le(1-\int_{\mathbb{B}_R}e^{8\pi \varphi}dx+o_\epsilon(R)\ri),
  \ena
  where $o_\epsilon(R)\ra 0$ as $\epsilon\ra 0$ for any fixed $R>0$.
  Letting $\epsilon\ra 0$ first, then $R\ra +\infty$, we have $I_1\ra
  0$. Secondly there holds
  \bna
  |I_2|\leq\sup_\mathbb{B}|\phi|\f{c_\epsilon}{\lambda_\epsilon}\int_\mathbb{B}
  u_\epsilon e^{(4\pi-\epsilon)u_{\epsilon,\beta}^2}dx.
  \ena
  It follows from Lemma \ref{L-5} and (\ref{t-z}) that $I_2\ra
  0$ as $\epsilon\ra 0$ first and then $R\ra+\infty$.

  Finally we can easily see that
  \bna
  I_3=\phi(\xi)
  \le(\int_{\mathbb{B}_R}e^{8\pi\varphi}dx+o_\epsilon(R)\ri)
  \ena
  for some $\xi\in \mathbb{B}_{Rr_\epsilon}$. Letting $\epsilon\ra 0$ first, then
  $R\ra +\infty$, we have $I_3\ra \phi(0)$.
  Combining all the above estimates, we finish the proof of
  the lemma. $\hfill\Box$\\

  For simplicity we denote
  $$\mathscr{L}_\alpha=-\Delta -\f{1}{(1-|x|^2)^2}-\alpha .$$
  Then we have the following:
  \begin{lemma}\label{L-9}
  The function sequence $c_\epsilon u_\epsilon$ converges to $G$ weakly
  in $W^{1,p}_{\rm loc}(\mathbb{B})$ for any $p\in (1,2)$,  strongly in $L^q(\mathbb{B})$ for any $q\geq 1$,
  and in $C^0(\overline{\mathbb{B}_r^c})$ for any $r\in (0,1)$, where $G$ is a Green function satisfying
  $\mathscr{L}_\alpha G=\delta_0$,
  where $\delta_0$ is the usual Dirac measure centered at $0\in \mathbb{B}$.
  \end{lemma}

 \noindent {\it Proof.} Note that
  $$\mathscr{L}_\alpha (c_\epsilon u_\epsilon)=f_\epsilon
  =\f{1}{\lambda_\epsilon}c_\epsilon u_\epsilon e^{(4\pi-\epsilon)u_\epsilon^2}.$$
  Let $v_\epsilon$ be a solution to
  \be\label{veps}\le\{\begin{array}{lll}
  \mathscr{L}_\alpha v_\epsilon=f_\epsilon\quad
  {\rm in}\quad\mathbb{B}_{1/2}\\
  [1.5ex] v_\epsilon=0\quad{\rm on}\quad \p\mathbb{B}_{1/2}
  \end{array}\ri.\ee
  By Lemma \ref{L-7}, $f_\epsilon$ is bounded in $L^1(\mathbb{B})$. By a result of Struwe \cite{Struwe1}, for any $q$, $1<q<2$, there holds
  \be\label{w1q}\|\nabla v_\epsilon\|_q\leq C\|f_\epsilon\|_1,\ee
  and there exists some $v_0\in W_0^{1,q}(\mathbb{B}_{1/2})$ such that
  \be\label{wq}v_\epsilon\rightharpoonup v_0\quad{\rm weakly\,\,in}\quad W_0^{1,q}(\mathbb{B}_{1/2}).\ee
  Take a cut-off function $\phi\in C_0^\infty(\mathbb{B})$ satisfying $0\leq \phi\leq 1$,
  $\phi\equiv 1$ on $\mathbb{B}_{1/8}$ and $\phi\equiv 0$ outside $\mathbb{B}_{1/4}$.
  Set $w_\epsilon=c_\epsilon u_\epsilon-\phi v_\epsilon$.
  It follows that
  $$\mathscr{L}_\alpha w_\epsilon=(1-\phi)f_\epsilon
  +\Delta \phi v_\epsilon+2\nabla \phi\nabla v_\epsilon.$$
  By (\ref{t-z}) and Lemma \ref{L-2}, $f_\epsilon$ is uniformly bounded in $\mathbb{B}\setminus\mathbb{B}_{1/16}$. While
  (\ref{w1q}) and the Sobolev embedding theorem imply that $v_\epsilon$ is bounded in $L^2(\mathbb{B}_{1/2})$.
  Then applying elliptic estimates to (\ref{veps}), we conclude that $v_\epsilon$ is bounded in
  $W^{2,2}(\mathbb{B}_{1/4}\setminus
  \mathbb{B}_{1/8})$, and thus $\nabla\phi\nabla v_\epsilon$ is bounded in $L^2(\mathbb{B})$. Therefore
  $\mathscr{L}_\alpha w_\epsilon$ is bounded in $L^2(\mathbb{B})$. Recalling Lemma \ref{L-2}, we have
  \bna
  \|w_\epsilon\|_{1,\alpha}^2=\langle w_\epsilon,\mathscr{L}_\alpha w_\epsilon\rangle_{L^2}
  \leq C\|w_\epsilon\|_{1,\alpha}\|\mathscr{L}_\alpha w_\epsilon\|_2.
  \ena
  This implies that $w_\epsilon$ is bounded in $\mathscr{H}$ and there exists some $w_0\in \mathscr{H}$ such that
  \be\label{wea}
  w_\epsilon\rightharpoonup w_0 \quad{\rm weakly\,\,in}\quad \mathscr{H}.
  \ee
  Let $G=\phi v_0+w_0$. Here we extend $v_0$ to be zero in $\mathbb{B}\setminus\mathbb{B}_{1/2}$. It follows from (\ref{wq}) and
  Lemma \ref{L-2} that
   $c_\epsilon u_\epsilon\ra G$ in $L^p(\mathbb{B})$ for any $p\geq 1$ and in $C^0(\mathbb{B}\setminus\mathbb{B}_r)$ for any
   $r>0$. Moreover we have
  $$\lim_{\epsilon\ra 0}\int_{\mathbb{B}} c_\epsilon u_\epsilon\mathscr{L}_\alpha \varphi dx=
   \int_{\mathbb{B}} G\mathscr{L}_\alpha \varphi dx,\quad\forall \varphi\in C_0^\infty(\mathbb{B}).$$
   This together with Lemma \ref{L-7} finishes the proof of the lemma. $\hfill\Box$\\

   Before ending this subsection, we decompose the Green function $G$. Since
   $$-\Delta \le(G+\f{1}{2\pi}\log r\ri)=\f{G}{(1-|x|^2)^2}+\alpha G\in L^p_{\rm loc}(\mathbb{B}),$$
   there holds
   \be\label{green}G=-\f{1}{2\pi}\log r+A_0+\widetilde{\psi},\ee
   where $\widetilde{\psi}\in C^1_{\rm loc}(\mathbb{B})$.

 \subsection{Neck analysis and upper bound estimate\\}

 In this subsection, we use the capacity estimate due to Y. Li \cite{Lijpde} to derive an upper bound of
 the supremum in (\ref{htm}). While in \cite{Wang-Ye}, this was done by G. Wang and D. Ye by using a result of Carleson-Chang
 \cite{CC}, which was employed originally by Li-Liu-Yang \cite{LLY} when deriving an upper bound of certain
 Trudinger-Moser functional for vector bundles on a compact Riemannian surface.

 \begin{lemma}\label{L-10} For any $r$, $0<r<1$, there holds
   $$\int_{\mathbb{B}_r}|\nabla u_\epsilon|^2dx=1+\f{1}{c_\epsilon^2}
   \le(\f{1}{2\pi}\log{r}-A_0+o_r(1)+o_\epsilon(1)\ri),$$
   where $o_\epsilon(1)\ra 0$ as $\epsilon\ra 0$, $o_r(1)\ra 0$ as $r\ra 0$.
   \end{lemma}

   \noindent{\it Proof}. In view of the Euler-Lagrange equation (\ref{Euler-Lagrange}), we have
   by using the divergence theorem
   \bna
   \int_{\mathbb{B}_r}|\nabla u_\epsilon|^2dx&=&-\int_{\mathbb{B}_r}u_\epsilon\Delta u_\epsilon dx
   +\int_{\p\mathbb{B}_r}u_\epsilon \f{\p u_\epsilon}{\p \nu}ds\\
   &=&\int_{\mathbb{B}_r}u_\epsilon\mathscr{L}_\alpha u_\epsilon dx+\int_{\p\mathbb{B}_r}u_\epsilon
   \f{\p u_\epsilon}{\p \nu} d s+\int_{\mathbb{B}_r}\f{u_\epsilon^2}{(1-|x|^2)^2}dx
   +\alpha \int_{\mathbb{B}_r}u_\epsilon^2dx\\
   &=&\int_{\mathbb{B}_r}\f{u_\epsilon^2}{\lambda_\epsilon}e^{(4\pi-\epsilon)u_\epsilon^2}dx
   +\int_{\p\mathbb{B}_r}u_\epsilon
   \f{\p u_\epsilon}{\p \nu} d s+\int_{\mathbb{B}_r}\f{u_\epsilon^2}{(1-|x|^2)^2}dx
   +\alpha \int_{\mathbb{B}_r}u_\epsilon^2dx.
   \ena
   Now we estimate the above four integrals respectively. It follows from Lemma \ref{L-9} and (\ref{t-z}) that
   $$\int_{\mathbb{B}_r}\f{u_\epsilon^2}{\lambda_\epsilon}e^{(4\pi-\epsilon)u_\epsilon^2}dx=
   1-\f{1}{c_\epsilon^2}\int_{\mathbb{B}\setminus\mathbb{B}_r}\f{(c_\epsilon u_\epsilon)^2}{\lambda_\epsilon}
   e^{(4\pi-\epsilon)u_\epsilon^2}dx=1-\f{1}{c_\epsilon^2}o_\epsilon(1).$$
   Moreover, Lemma \ref{L-9} and (\ref{green}) lead to
   $$\int_{\p\mathbb{B}_r}u_\epsilon
   \f{\p u_\epsilon}{\p \nu} d s=\f{1}{c_\epsilon^2}\le(\int_{\p\mathbb{B}_r}G\f{\p G}{\p r}ds+o_\epsilon(1)\ri)
   =\f{1}{c_\epsilon^2}\le(\f{1}{2\pi}\log r-A_0+o_r(1)\ri),$$
   $$\int_{\mathbb{B}_r}\f{u_\epsilon^2}{(1-|x|^2)^2}dx=\f{1}{c_\epsilon^2}\le(
   \int_{\mathbb{B}_r}\f{G^2}{(1-|x|^2)^2}dx+o_\epsilon(1)\ri)=\f{o_r(1)+o_\epsilon(1)}{c_\epsilon^2},$$
   and
   $$\int_{\mathbb{B}_r}u_\epsilon^2dx=\f{1}{c_\epsilon^2}\le(\int_{\mathbb{B}_r}G^2dx+o_\epsilon(1)\ri)=
   \f{o_r(1)+o_\epsilon(1)}{c_\epsilon^2}.$$
   Combining all the above estimates, we finish the proof of the lemma. $\hfill\Box$

   \begin{lemma}\label{L-11}
   For two positive numbers $\delta$ and $R$ with $\delta>Rr_\epsilon$, there holds
   $$\int_{\mathbb{B}_\delta\setminus\mathbb{B}_{Rr_\epsilon}}|\nabla u_\epsilon|^2dx=1+
   \f{1}{c_\epsilon^2}\le(-\f{\log R}{2\pi}+\f{\log\delta}{2\pi}-\f{\log\pi}{4\pi}+\f{1}{4\pi}-A_0+o_\delta(1)+O(\f{1}{R^2})+o_\epsilon(1)\ri).$$
   \end{lemma}

   \noindent{\it Proof}. By (\ref{phibubble}) and (\ref{v-r}), we have
   \bna
   \int_{\mathbb{B}_{Rr_\epsilon}}|\nabla u_\epsilon|^2dx&=&\f{1}{c_\epsilon^2}\le(\int_{\mathbb{B}_{R}}|\nabla\varphi(y)|^2dy+
   o_\epsilon(1)\ri)\\
   &=&\f{1}{c_\epsilon^2}\le(\f{\log R}{2\pi}+\f{\log\pi}{4\pi}-\f{1}{4\pi}+O(\f{1}{R^2})+o_\epsilon(1)\ri).
   \ena
   This together with Lemma \ref{L-10} implies the lemma. $\hfill\Box$\\

   Let $0<s<r<1$ and $a,b\in\mathbb{R}$. The function $h:\mathbb{B}_r\setminus\mathbb{B}_s\ra \mathbb{R}$ defined by
   $$
   h(x)=\f{b\log\f{|x|}{s}+a\log\f{r}{|x|}}{\log\f{r}{s}},
   $$
   is harmonic on the planar domain  $\mathbb{B}_r\setminus\mathbb{B}_s$. Obviously $h$ has boundary values
   $$\label{bdv}h|_{\p\mathbb{B}_s}=a,\,\,\,h|_{\p \mathbb{B}_r}=b. $$
   Moreover we have
   \be\label{enr-h}\int_{\mathbb{B}_r\setminus\mathbb{B}_s}|\nabla h|^2dx=\f{2\pi(b-a)^2}{\log\f{r}{s}}.\ee
   Define a function space associated with $h$ as
   $$\mathcal{W}=\mathcal{W}(h,r,s)=\le\{u\in W^{1,2}(\mathbb{B}_r\setminus\mathbb{B}_s)\,|\,
   u-h\in W_0^{1,2}(\mathbb{B}_r\setminus\mathbb{B}_s)\ri\}.$$
   By a variational direct method, one can see that the infimum
   $$\inf_{u\in\mathcal{W}}\int_{\mathbb{B}_r\setminus\mathbb{B}_s}|\nabla u|^2dx$$
   can be attained by the above harmonic function $h$. In fact we have proved the following:

   \begin{lemma}\label{L-12}
   Let $0<s<r<1$, $a,b\in\mathbb{R}$, and $h$, $\mathcal{W}$ be given as above.
   There holds
   $$\inf_{u\in\mathcal{W}}\int_{\mathbb{B}_r\setminus\mathbb{B}_s}|\nabla u|^2dx
   =\f{2\pi(b-a)^2}{\log\f{r}{s}}.$$
   \end{lemma}

   This lemma can be used to derive the following:

   \begin{lemma}\label{L-13} Assume $0<\delta<1$, $R>0$ and $\epsilon$ is sufficiently small. Then there holds
   $$\int_{\mathbb{B}_\delta\setminus\mathbb{B}_{Rr_\epsilon}}
   |\nabla u_\epsilon|^2dx\geq
   \f{2\pi(b_\epsilon-a_\epsilon)^2}{\log\f{\delta}{Rr_\epsilon}},$$
   where $a_\epsilon$ and $b_\epsilon$ are defined as
   \bna
   a_\epsilon&=&u_\epsilon|_{\p \mathbb{B}_{Rr_\epsilon}}=c_\epsilon+\f{1}{c_\epsilon}\le(
   -\f{1}{2\pi}\log R-\f{1}{4\pi}\log\pi+O(\f{1}{R^2})+o_\epsilon(1) \ri),\\
   b_\epsilon&=&u_\epsilon|_{\p \mathbb{B}_{\delta}}=\f{1}{c_\epsilon}\le(-\f{1}{2\pi}\log\delta+A_0+o_\delta(1)+o_\epsilon(1)\ri).
   \ena
   \end{lemma}

   \noindent{\it Proof}. Substitute $a_\epsilon$, $b_\epsilon$, $Rr_\epsilon$ and $\delta$ for
   $a$, $b$, $s$ and $r$ respectively in Lemma \ref{L-12}. Let
   $$h_\epsilon(x)=\f{b_\epsilon\log\f{|x|}{Rr_\epsilon}+a_\epsilon\log\f{\delta}{|x|}}{\log\f{\delta}{Rr_\epsilon}}.$$
   Then $h_\epsilon|_{\p \mathbb{B}_{Rr_\epsilon}}=u_\epsilon|_{\p \mathbb{B}_{Rr_\epsilon}}$ and
   $h_\epsilon|_{\p \mathbb{B}_{\delta}}=u_\epsilon|_{\p \mathbb{B}_{\delta}}$. Hence we have $u_\epsilon-h_\epsilon\in
   W_0^{1,2}(\mathbb{B}_{\delta}\setminus\mathbb{B}_{Rr_\epsilon})$ and
   $$\int_{\mathbb{B}_\delta\setminus\mathbb{B}_{Rr_\epsilon}}
   |\nabla u_\epsilon|^2dx\geq\inf_{v\in\mathcal{W}(h_\epsilon,\delta,Rr_\epsilon)}\int_{\mathbb{B}_\delta\setminus\mathbb{B}_{Rr_\epsilon}}
   |\nabla v|^2dx=\int_{\mathbb{B}_\delta\setminus\mathbb{B}_{Rr_\epsilon}}
   |\nabla h_\epsilon|^2dx.$$
   This together with an obvious analog of (\ref{enr-h}) concludes the lemma.$\hfill\Box$ \\

   A straightforward calculation shows
   \bea\nonumber
   2\pi(b_\epsilon-a_\epsilon)^2&=&2\pi \le\{c_\epsilon+\f{1}{c_\epsilon}
   \le(-\f{1}{2\pi}\log R+\f{1}{2\pi}\log\delta-\f{1}{4\pi}\log\pi-A_0+o(1)\ri)\ri\}^2\\
   &=&2\pi c_\epsilon^2\le\{1+\f{1}{c_\epsilon^2}\le(-\f{1}{\pi}\log R+\f{1}{\pi}\log\delta
   -\f{1}{2\pi}\log\pi-2A_0+o(1)\ri)\ri\}.\label{b-a}
   \eea
   Here and in the sequel $o(1)\ra 0$ as $\epsilon\ra 0$ first, then $R\ra+\infty$, and finally $\delta\ra 0$. Also we have
   \be\label{log}\log\f{\delta}{Rr_\epsilon}=\log\delta-\log R-\log\f{\sqrt{\lambda_\epsilon}}{c_\epsilon}
   +(2\pi-\epsilon/2)c_\epsilon^2.\ee
   Combining Lemma \ref{L-11}, Lemma \ref{L-13}, (\ref{b-a}) and (\ref{log}), we obtain
   \bna
   &&1+
   \f{1}{c_\epsilon^2}\le(-\f{\log R}{2\pi}+\f{\log\delta}{2\pi}-\f{\log\pi}{4\pi}+\f{1}{4\pi}-A_0+o(1)\ri)\\
   &&\geq
   \f{1+\f{1}{c_\epsilon^2}\le(-\f{\log R}{\pi}+\f{\log\delta}{\pi}
   -\f{\log\pi}{2\pi}-2A_0+o(1)\ri)}{1-\f{\epsilon}{4\pi}+\f{1}{c_\epsilon^2}
   \le(-\f{\log R}{2\pi}+\f{\log\delta}{2\pi}-\f{1}{2\pi}\log\f{\sqrt{\lambda_\epsilon}}{c_\epsilon}\ri)}.
   \ena
   This leads to
   \bna
   &&1+\f{1}{c_\epsilon^2}\le(-\f{\log R}{\pi}+\f{\log \delta}{\pi}-\f{\log\pi}{4\pi}+\f{1}{4\pi}-A_0
   -\f{1}{2\pi}\log\f{\sqrt{\lambda_\epsilon}}{c_\epsilon}+o(1)\ri)\\
   &&\geq 1+\f{1}{c_\epsilon^2}\le(-\f{\log R}{\pi}+\f{\log\delta}{\pi}
   -\f{\log\pi}{2\pi}-2A_0+o(1)\ri).
   \ena
   It then follows that
   $$\f{1}{2\pi}\log\f{\sqrt{\lambda_\epsilon}}{c_\epsilon}\leq \f{\log \pi}{4\pi}+\f{1}{4\pi}+A_0+o(1).$$
   Therefore
   $$\label{u-p}\f{\lambda_\epsilon}{c_\epsilon^2}\leq \pi e^{1+4\pi A_0+o(1)}.$$
   This together with Lemma \ref{L-6} implies the following:

   \begin{proposition}\label{P-13} Under the assumption that $c_\epsilon=\max_{\mathbb{B}}
   u_\epsilon\ra +\infty$ as $\epsilon\ra 0$, there holds
   \be\label{pro}\sup_{u\in\mathscr{H},\,\|u\|_{1,\alpha}\leq 1}\int_{\mathbb{B}}e^{4\pi u^2}dx=
   \lim_{\epsilon\ra 0} \int_\mathbb{B}e^{(4\pi-\epsilon)u_\epsilon^2}dx\leq\pi+\pi e^{1+4\pi A_0}.\ee
   \end{proposition}

   \subsection{Test function computation\\}

     In this subsection, we construct a sequence of test functions $\phi_\epsilon\in\mathscr{H}$ such that
     $\|\phi_\epsilon\|_{1,\alpha}\leq 1$ and if $\epsilon$ is chosen sufficiently small, there holds
     \be\label{int}\int_{\mathbb{B}}e^{4\pi \phi_\epsilon^2}dx>\pi+\pi e^{1+4\pi A_0}.\ee
     By Proposition \ref{P-13}, this would contradicts (\ref{pro}) unless $c_\epsilon$ is bounded. Therefore we get the desired extremal function
     and complete the proof of Theorem \ref{thm1}.

       We set
     \be\label{ppp}\phi_\epsilon(x)=\le\{
     \begin{array}{llll}
     &c+\f{-\f{1}{4\pi}\log(1+\pi\f{|x|^2}{\epsilon^2})+B}{c}\quad &{\rm for} & |x|\leq R\epsilon\\[1.5ex]
     &\f{G(x)}{c}\quad &{\rm for} &  R\epsilon<|x|\leq 1,
     \end{array}
     \ri.\ee
     where $R=-\log\epsilon$, $B$ and $c$ are constants to be determined
     later.
     We now require
     \be\label{cont}
     c+\f{1}{c}\le(-\f{1}{4\pi}\log(1+\pi R^2)+B\ri)=\f{1}{c}G\mid_{\p\mathbb{B}_{R\epsilon}}
     =\f{1}{c}\le(-\f{1}{2\pi}\log (R\epsilon)+A_{0}+O(R\epsilon)\ri),
     \ee
     which gives
     \be\label{2pic2-1}
     2\pi c^2=-\log\epsilon-2\pi B+2\pi A_{0}+\f{1}{2}\log \pi
     +O(\f{1}{R^2}).
     \ee
     Clearly, (\ref{ppp}) and (\ref{cont}) imply that $\phi_\epsilon\in W^{1,2}_{\rm loc}(\mathbb{B})$. While
     in view of (\ref{wea}),
     $G$ coincides with $w_0\in\mathscr{H}$ on $\mathbb{B}\setminus\mathbb{B}_{{1}/{2}}$. Hence
     $\phi_\epsilon-{w_0}/{c}\in W_0^{1,2}(\mathbb{B})$, which immediately leads to the fact that $\phi_\epsilon\in\mathscr{H}$.

     Since $\phi_\epsilon\in \mathscr{H}$, we have by integration by parts, Lemma \ref{L-9} and (\ref{green}) that
     \bna
     \int_{\mathbb{B}\setminus\mathbb{B}_{R\epsilon}}\le(|\nabla \phi_\epsilon|^2-\f{\phi_\epsilon^2}
     {(1-|x|^2)^2}-\alpha \phi_\epsilon^2\ri)dx&=&\f{1}{c^2}\int_{\p\mathbb{B}_{R\epsilon}}G\f{\p G}{\p\nu}ds+
     \f{1}{c^2}\int_{\mathbb{B}\setminus\mathbb{B}_{R\epsilon}}G\mathscr{L}_\alpha Gdx\\
     &=&\f{1}{c^2}\le(-\f{1}{2\pi}\log(R\epsilon)+A_0+O(\f{1}{R^2})\ri).
     \ena
     Also a straightforward calculation gives
     \bna
     \int_{\mathbb{B}_{R\epsilon}}|\nabla\phi_\epsilon|^2dx=\f{1}{c^2}
     \le(\f{\log R}{2\pi}+\f{\log\pi}{4\pi}-\f{1}{4\pi}+O(\f{1}{R^2})\ri).
     \ena
     Hence
     \bna
     \|\phi_\epsilon\|_{1,\alpha}^2&\leq&\int_{\mathbb{B}\setminus\mathbb{B}_{R\epsilon}}\le(|\nabla \phi_\epsilon|^2-\f{\phi_\epsilon^2}
     {(1-|x|^2)^2}-\alpha \phi_\epsilon^2\ri)dx+\int_{\mathbb{B}_{R\epsilon}}|\nabla\phi_\epsilon|^2dx\\
     &=&\f{1}{c^2}\le(-\f{1}{2\pi}\log\epsilon+A_0+\f{\log\pi}{4\pi}-\f{1}{4\pi}+O(\f{1}{R^2})\ri).
     \ena
     We set
     $$\int_{\mathbb{B}\setminus\mathbb{B}_{R\epsilon}}\le(|\nabla \phi_\epsilon|^2-\f{\phi_\epsilon^2}
     {(1-|x|^2)^2}-\alpha \phi_\epsilon^2\ri)dx+\int_{\mathbb{B}_{R\epsilon}}|\nabla\phi_\epsilon|^2dx=1,$$
     which implies $\|\phi_\epsilon\|_{1,\alpha}\leq 1$ and
     \be\label{c}2\pi c^2=-\log\epsilon+2\pi A_0+\f{1}{2}\log\pi-\f{1}{2}+O(\f{1}{R^2}).\ee
     Combining (\ref{2pic2-1}) and (\ref{c}), we obtain
     \be\label{B}B=\f{1}{4\pi}+O(\f{1}{R^2}).\ee

     We now derive the estimate (\ref{int}). It is clear that
     \bna
     \int_{\mathbb{B}\setminus\mathbb{B}_{R\epsilon}}e^{4\pi \phi_\epsilon^2}dx&\geq&
     \int_{\mathbb{B}\setminus\mathbb{B}_{R\epsilon}}(1+4\pi\phi_\epsilon^2)dx\\
     &=&\pi+\f{4\pi}{c^2}\le(\int_\mathbb{B}G^2dx+O(\f{1}{R^2})\ri).
     \ena
     By (\ref{c}) and (\ref{B}), there holds on $\mathbb{B}_{R\epsilon}$,
     \bna\phi_\epsilon^2&\geq&c^2+2B-\f{1}{2\pi}\log\le(1+\pi\f{|x|^2}{\epsilon^2}\ri)\\
     &=&-\f{1}{2\pi}\log\le(1+\pi\f{|x|^2}{\epsilon^2}\ri)-\f{1}{2\pi}\log\epsilon+A_0+\f{\log\pi}{4\pi}+
     \f{1}{4\pi}+O(\f{1}{R^2}).
     \ena
     This leads to
     \bna
     \int_{\mathbb{B}_{R\epsilon}}e^{4\pi \phi_\epsilon^2}dx
     &\geq&e^{1+4\pi A_0+\log\pi+O(\f{1}{R^2})}\int_{\mathbb{B}_R}\f{1}
     {(1+\pi|x|^2)^2}dx\\
     &=& \pi e^{1+4\pi A_0}\le(1+O(\f{1}{R^2})\ri).
     \ena
     Since $1/R^2=o(1/c^2)$, we obtain
     \bna
     \int_{\mathbb{B}}e^{4\pi \phi_\epsilon^2}dx&=&\int_{\mathbb{B}\setminus\mathbb{B}_{R\epsilon}}e^{4\pi \phi_\epsilon^2}dx
     +\int_{\mathbb{B}_{R\epsilon}}e^{4\pi \phi_\epsilon^2}dx\\
     &\geq&\pi+\pi e^{1+4\pi A_0}+\f{4\pi}{c^2}\le(\int_{\mathbb{B}}G^2dx+o(1)\ri).
     \ena
     This gives the desired estimate (\ref{int}) provided that $\epsilon$ is sufficiently small. \\

    {\bf Acknowledgements}. Y. Yang is supported by the National Science Foundation of China (Grant No.11171347 and Grant
 No. 11471014). X. Zhu is supported by the National Science Foundation of China (Grant No. 41275063 and Grant
 No. 1140575).


\end{document}